\documentclass[12pt,reqno]{amsart}
\usepackage{amssymb,delarray}
\usepackage{graphicx}

\usepackage{epsf,amsfonts,amscd,latexsym}

\def\leq{\leqslant}

\newtheorem{thm}{Theorem}[section]

\newtheorem{cor}[thm]{Corollary}

\begin{document}
\baselineskip=14.pt plus 2pt 

\title[
]
{On the fundamental groups of the complements of Hurwitz curves}
\author[Olga Kulikova]{O.V.~Kulikova}
\address{Moscow State University
} \email{olga.kulikova@mail.ru}

\dedicatory{} \subjclass{}
\begin{abstract}
It is proved that the commutator subgroup of the fundamental group
of the complement of any plane affine irreducible Hurwitz curve
(respectively, any plane affine irreducible pseudoholomorphic
curve) is finitely presented. It is shown that there exists a
pseudo-holomorphic curve (a Hurwitz curve) in $\mathbb C\mathbb
P^2$ whose fundamental group of the complement is not Hopfian and,
respectively, this group is not residually finite. In addition, it
is proved that there exist an irreducible nonsingular algebraic
curve $C\subset \mathbb C^2$ and a bi-disk $D\subset \mathbb C^2$
such that the fundamental group $\pi_1(D\setminus C)$ is not
Hopfian.\end{abstract}

\maketitle
\setcounter{tocdepth}{2}


\def\st{{\sf st}}

\setcounter{section}{-1}
\section{Introduction}

The notion of Hurwitz curves with respect to a linear projection
of the projective plane $\mathbb C\mathbb P^2$ to $\mathbb
C\mathbb P^1$ was introduced in \cite{Moi2}  and is a natural
generalization of the notion of the plane algebraic curves (in
\cite{Moi2}, Hurwitz curves are called "semi-algebraic curves"). A
precise definition of Hurwitz curves can be found, for example, in
\cite{Kh-Ku}. Roughly speaking, Hurwitz curves in  $\mathbb
C\mathbb P^2$ imitate the behavior of plane algebraic curves with
respect to the pencil of complex lines defining the projection. In
particular, they look like analytic curves in neighborhoods of
critical points of the projection.

Hurwitz curves play an important role in symplectic geometry. In
particular, Auroux and Katzarkov (see \cite{Au}, \cite{Au-Ka})
proved that a compact symplectic 4-manifold $(X,\omega)$ with
symplectic form $\omega$, whose class $[\omega]\in H^2(X,\mathbb
Z)$ and for which an $\omega$-compatible almost complex structure
$J$ is chosen, can be presented as an approximately holomorphic
generic covering $f_k:X\to \mathbb C\mathbb P^2$, $k>>0$, branched
over a cuspidal Hurwitz curve $\bar H_k$ (maybe with negative
nodes), where $f_k$ is given by three sections of the line bundle
$L^{\otimes k}$ and $L$ is a line bundle on $X$ whose first Chern
class is $[\omega]$. Therefore the investigation of the
fundamental groups $\pi_1(\mathbb C\mathbb P^2\setminus \bar H)$
of the complements of Hurwitz curves $\bar H$ is very important
for symplectic geometry.

In \cite{Ku1}, a class $\mathcal C$ of groups, called $C$-groups,
was defined. By definition, this class consists of the groups $G$
which are given by finite presentations of the following form: for
some integer $m$, a subset
$$J=\{ \, (i,j,k)\in \mathbb Z^3 \,
\mid \, 1\leq i,j\leq m,\, \, 1\leq k\leq h(i,j)\} ,$$ where
$h:\{1,\dots , m\}^2 \to \mathbb Z$ is a function, and a subset
$W=\{ w_{i,j,k}\in \mathbb F_m\, \mid \, (i,j,k)\in J\} $ of words
in a free group $\mathbb F_m$ generated by an alphabet
$\{x_1,\dots, x_m\}$   (it is possible that
$w_{i_1,j_1,k_1}=w_{i_2,j_2,k_2}$ for $(i_1,j_1,k_1)\neq
(i_2,j_2,k_2)$), a group $G\in \mathcal C$ possesses the
presentation
\begin{equation} \label{zero}
G_W=<x_1,\dots ,x_m \, \mid \, x_i= w_{i,j,k}^{-1} x_jw_{i,j,k} ,
\, \, w_{i,j,k}\in W\, >.
\end{equation}
A $C$-group $G$ is called {\it irreducible} if $G/G'\simeq \mathbb
Z$, where $G'=[G,G]$ is the commutator subgroup.

Denote by $\varphi_W: \mathbb F_m\to G_W$  the canonical
epimorphism. The elements $\varphi_{W}(x_i)\in G_W$, $1\leq i\leq
m$, and the elements conjugated to them are called the {\it
$C$-generators} of the $C$-group $G=G_W$. Let $f:G_1\to G_2$ be a
homomorphism of $C$-groups. It is called a {\it $C$-homomorphism}
if the images of the $C$-generators of $G_1$ under $f$ are
$C$-generators of the $C$-group $G_2$.  We will distinguish
$C$-groups up to $C$-isomorphisms.

Note that the class $\mathcal C$ contains the subclasses $\mathcal
K$ and $\mathcal L$, respectively, of the knot and link groups
given by Wirtinger presentation. In \cite{Ku1}, it was proved that
the class $\mathcal C$ coincides with the class of the fundamental
groups of the complements of orientable closed surfaces in the
4-dimensional sphere $S^4$ (with generalized Wirtinger
presentation). Besides, it follows from Theorem 1.14 in \cite{Ku}
and Theorem 2.1 in \cite{Kh-Ku} that for any $C$-group $G$ there
are an affine nonsingular algebraic curve $C\subset \mathbb C^2$
and a bi-disk $D=\{ \, (z,w)\in \mathbb C^2\, \, \mid \, \, \mid
z\mid \leq 1,\, \, \mid w\mid \leq 1\, \}$ such that $G\simeq
\pi_1(D\setminus C)$.

Let $H\subset \mathbb C^2=\mathbb C\mathbb P^2\setminus
L_{\infty}$ be an affine Hurwitz curve of degree $m$, that is,
$H=\bar H\cap \mathbb C^2$, where $\bar H$ is a Hurwitz curve in
$\mathbb C\mathbb P^2$ with respect to some pencil of lines and a
line $L_{\infty}$ is a member of the pencil being in general
position with respect to $\bar H$, and the intersection number
$\bar H\cdot L_{\infty}=m$. Denote by $\mathcal H=\{ \,
\pi_1(\mathbb C^2\setminus H)\, \}$ the class of the fundamental
groups of the complements of the affine Hurwitz curves $H$. If
$\bar H$ is a Hurwitz curve  of degree $m$, then the Zariski --
van Kampen presentation of $\pi_1(\mathbb C^2\setminus H)$
(defined by the pencil of lines) is a presentation of the form
(\ref{zero}) such that the words $w_{i,i,1}=x_1\dots x_m$,
$i=1,\dots,m$, belong to $W$, i.e., the element $x_1\dots x_m$
belongs to the center of $G_W$. In \cite{Ku}, it was proved that
if a $C$-group $G_W$ is given by presentation (\ref{zero}) such
that the element $x_1\dots x_m$ belongs to the center of $G_W$,
then for some $k\in \mathbb N$ there is a Hurwitz curve $\bar
H\subset \mathbb C\mathbb P^2$ of degree $M=m2^k$ such that
$G_W=\pi_1(\mathbb C^2\setminus H)$ and the element $(x_1\dots
x_m)^{2^k}\in G_W$ corresponds to a circuit around $L_{\infty}$.
One can add the generators $x_{m+1},\dots , x_{2^km}$ and
relations $x_i=x_{i+m}$ for $i=1,\dots ,(2^k-1)m$ to presentation
(\ref{zero}) of $G_W$ and obtain a group isomorphic to $G_W$. Note
that relations $x_i=x_j$ are $C$-relations, since they can be
written as $x_j=x_j^{-1}x_ix_j$. Therefore the class $\mathcal H$
coincides with the subclass of $\mathcal C$ consisting of the
groups $G$ which possess presentations (\ref{zero}) for some $m$
and such that the elements $x_1\dots x_m$ (for some $m$ for each
group $G\in \mathcal H$) belong to the centers of these groups,
since we consider $C$-groups up to $C$-isomorphisms. A group $G\in
\mathcal H$ will be called {\it a Hurwitz $C$-group of degree $m$}
if $C$-generators $x_1,\dots ,x_m$ of $G$ generate the group $G$
and the element $x_1\dots x_m$ belongs to the center of $G$ (note
that the degree of a Hurwitz $C$-group $G$ is not defined
canonically and depends on its $C$-presentation).

As is known the commutator subgroups $G'$ of a lot of irreducible
$C$-groups $G$ are not finitely generated (in particular, $G'$ of
$G\in \mathcal K$ is finitely generated group iff $G$ is the group
of a fibred knot (see \cite{S}), and, moreover, in the case of a
fibred knot $G'$ is a free group). One of the main results of this
article is the following theorem.

\begin{thm} \label{main} Let $G=\pi_1(\mathbb C^2\setminus H)\in \mathcal H$, where $\bar H$ is an
irreducible Hurwitz curve. Then the commutator subgroup $G'=[G,G]$
is a finitely presented group.
\end{thm}

Note that Theorem \ref{main} is a generalization of the similar
result in algebraic case (see \cite{Ku2}).

For any $C$-group $G=G_W$ with presentation (\ref{zero}), denote
by $\nu :G\to \mathbb F_1$ the natural $C$-homomorphism and
$N=\ker \nu$. Since $N=[G,G]$ for an irreducible $C$-group $G$,
Theorem \ref{main} is a simple consequence of the following

\begin{thm} \label{Main} For any Hurwitz $C$-group $G\in \mathcal H$ the group $N$ is finitely presented.
\end{thm}
In difference to \cite{Ku2}, the proof of Theorem \ref{Main} given
given in section 1 is purely algebraic.

Let $J$ be an almost complex structure in $\mathbb C\mathbb P^2$
compactible with Fubini -- Studi symplectic form and $\bar H$ be a
$J$-holomorphic curve in $\mathbb C\mathbb P^2$. If we chose a
pencil of pseudo-holomorphic lines, then, by Zariski -- van Kampen
Theorem, a presentation of $\pi_1(\mathbb C\mathbb P^2\setminus
(\bar H\cup L_{\infty}))$ is defined by braid monodromy
factorization of $\bar H$ with respect to the chosen pencil, where
$L_{\infty}$ is one of the $J$-lines of the pencil being in
general position with respect to $\bar H$. Therefore similar to
the case of Hurwitz curves, it is easy to show (see the proof of
Theorem 6.1 in \cite{Ku}) that $\pi_1(\mathbb C\mathbb
P^2\setminus (\bar H\cup L_{\infty}))$ is a Hurwitz $C$-group.
Thus,  we have
\begin{cor} \label{cor2} Let $\bar H$ be an irreducible pseudoholomorphic  curve in $\mathbb C\mathbb P^2$.
Then the commutator subgroup $G'$ of $G=\pi_1(\mathbb C\mathbb
P^2\setminus (\bar H\cup L_{\infty}))$ is a finitely presented
group.
\end{cor}

Let $\bar H$ be a Hurwitz curve of degree $m$. To obtain a
presentation of $\pi_1(\mathbb C\mathbb P^2\setminus \bar H)$ from
Zariski -- van Kampen presentation (\ref{zero}) of the group
$\pi_1(\mathbb C^2\setminus H)$, it is sufficient to add the
additional relation $x_1\dots x_m=1$ (the element $x_1\dots x_m\in
\pi_1(\mathbb C^2\setminus H)$ corresponds to a circuit around the
line $L_{\infty}$).

For any Hurwitz $C$-group $G=G_W$ of degree $m$ given by
presentation (\ref{zero}), denote by
$$\bar G_{m,k}=G_W/\{ (x_1\dots x_m)^k=1\}$$ and call $\bar G_{m,k}$
a {\it projective Hurwitz group of degree} $mk$. It is easy to see
that the homomorphism $\nu$ induces the homomorphism $\nu_{mk}
:\bar G_W\to \mathbb Z/mk\mathbb Z$. Put $\bar N_{m,k}=\ker
\nu_{mk}$. The following theorem is a particular case of Corollary
2.8. in \cite{M-K-S}.
\begin{thm} \label{Main2} For any projective Hurwitz group $\bar G_{m,k}$ of degree $mk$
the group $\bar N_{m,k}$ is finitely presented.
\end{thm}

Since for an irreducible Hurwitz curve (respectively, for an
irreducible pseudo-holomorphic curve) $\bar H$ of $\deg \bar H=m$
the commutator subgroup $G'$ of $G=\pi_1(\mathbb C\mathbb
P^2\setminus \bar H)$, given by Zariski -- van Kampen
presentation, coincides with $\bar N_{m,1}$, we have

\begin{cor} \label{cor} Let $\bar H\subset \mathbb C\mathbb P^2$
be an irreducible Hurwitz curve (respectively, pseudo-holomorphic
curve). Then the commutator subgroup $G'$ of $G=\pi_1(\mathbb
C\mathbb P^2\setminus \bar H)$ is a finitely presented group.
\end{cor}

 Let $C$ be a plane algebraic curve. In \cite{Zar},
O. Zariski formulated the following question:

{\it Is $G=\pi_1(\mathbb C\mathbb P^2 \setminus C)$ a residually
finite group?}
\newline It is natural to ask the same question in the local case, i.e.,
if $G=\pi_1(D \setminus C)$, where $D$ is a bi-disk in $\mathbb
C^2$, and in the cases of Hurwitz $C$-groups and projective
Hurwitz groups.

\begin{thm} \label{Zar} There are
\begin{itemize}
\item[$(i)$] an irreducible $C$-group $\tilde G$, \item[$(ii)$] a
Hurwitz $C$-group $G$ of degree three for which $G/G'\simeq
\mathbb Z^2$
\end{itemize}
such that $\tilde G$, $G$, and the projective Hurwitz groups $\bar
G_{3,k}$, $k\in \mathbb N$, associated with $G$, are non Hopfian.
In particular, they are not residually finite groups.
\end{thm}

Applying Theorems 1.14, 6.2 in \cite{Ku} and Theorem 2.1 in
\cite{Kh-Ku}, we obtain
\begin{cor} \label{co} There are
\begin{itemize}
\item[$(i)$] an irreducible nonsingular algebraic curve $C\subset
\mathbb C^2$ and a bi-disk $D=\{ \, (z,w)\in \mathbb C^2\, \, \mid
\, \, \mid z\mid \leq 1,\, \, \mid w\mid \leq 1\, \}$;
\item[$(ii)$] a Hurwitz curve $\bar H\subset \mathbb C\mathbb P^2$
consisting of two irreducible components
\end{itemize}
whose groups $\pi_1(D\setminus C)$, $\pi_1(\mathbb C\mathbb
P^2\setminus \bar H)$, and  $\pi_1(\mathbb C^2\setminus  H)$ are
non-Hopfian and, in particular, they are  not residually finite
groups.
\end{cor}
Note that, contrary to Corollary \ref{co} (i), if $C\subset
\mathbb C^2$ is a non-singular algebraic curve meeting
transversally the line at infinity, then, by Zariski Theorem,
$\pi_1(\mathbb C^2\setminus C)$ is an abelian group. Note also
that it follows from the proof of Theorem 6.2 in \cite{Ku} that
the Hurwitz curve $\bar H$ in Corollary \ref{co} (ii) can be
chosen such that all its singular points are simple triple points,
i.e., locally they can be given by equation $w(w^2-z^2)=0$.
Applying rescaling, one can assume that $\bar H$ is a symplectic
surface. Therefore we have
\begin{cor}  There exists a pseudo-holomorphic curve $\bar H\subset \mathbb C\mathbb P^2$
consisting of two irreducible components, having simple triple
points as its singularities and
 such that
 $\pi_1(\mathbb C\mathbb P^2\setminus \bar
H)$  is not Hopfian.
\end{cor}

The proof of Theorem \ref{Zar} is given in section 2.
\section{Proof of Theorem \ref{Main}}

Obviously, Theorem \ref{Main} follows from

\begin{thm} \label{Main3} Let a group $G$ given by a presentation
$$< x_1, \dots , x_n \, \, |\, \,  (x_1 \dots
x_n)x_j = x_j(x_1\dots x_n), \, \, j\in \{1,\dots ,n\}, R
>,$$
where $R$ is a set of words on the alphabet $X = \{ x_1, \dots ,
x_n \}$ such that there exists a homomorphism $\phi : G\to \mathbb
F_1$ mapping each $x_j$ to a generator $x$ of $\mathbb F_1$.
Denote by $N$ the kernel of $\phi$. If $R$ is finite, then $N$ is
finitely presented. In particular, if $R$ is empty, then $N$ is a
finitely generated free group.
\end{thm}
{\it Proof.} We have the exact sequence
$$ 1\to N\to G\to \mathbb F_1\to 1. $$
By the Tietze theorem, 
\begin{align*}
G\simeq & <  x_1, \dots , x_n, y \, \, |\,  \, , y = x_1\dots x_n,
(x_1\dots x_n) x_j =
x_j(x_1\dots x_n), \\
& j\in \{1,\dots ,n\},  R
>\simeq \\ & <  x_1, \dots , x_n, y \, |\,  y =
x_1\dots x_n, y x_j = x_j y, j\in \{1,\dots,n\},  R >\simeq \\
& \simeq <  x_2, \dots , x_n, y \, |\, y x_j = x_jy,\,  j\in
\{2,\dots ,n\}, \bar R >, \end{align*} where $\bar R$ is obtained
from $R$ if one changes $x_1$ by $y(x_2\dots x_n)^{-1}$ in every
word of $R$.

In terms of the last presentation, the homomorphism $\phi$ is
organized as follows: $\phi(x_j)=x$ for each $j\in\{2,\dots ,n\}$
and $\phi(y)=x^n$.

To find a finite presentation for $N$ let us use  Reidemeister --
Schreier method (see, for example, \S 2.3 \cite{M-K-S}). The
elements $x_n^k$ for $k\in \mathbb Z$ can be chosen as Schreier
representatives of cosets of $N$ in $G$. Then the group $N$ is
generated by
$$ a_{k,j} = x_n^kx_j\overline{x_n^k x_j}^{\, -1} =
x_n^kx_jx_n^{-(k+1)},$$ where  $j\in\{2,\dots ,n-1\}, k\in \mathbb
Z$, and the elements
$$ a_{k,n} = x_n^ky\overline{x_n^ky}^{\, -1} =
x_n^kyx_n^{-(k+n)},$$ where  $k\in \mathbb Z$. The relations
$$ yx_j = x_jy, \quad j=2,\dots ,n,$$
give rise to the relations
\begin{equation}\label{eq1}
a_{k,n}a_{k+n,j}a_{k+1,n}^{-1}a_{k,j}^{-1}=1
\end{equation}
for $
 j\in \{2,\dots ,n-1\}$ and $k\in \mathbb Z$ and the relations
 \begin{equation} \label{eq2}
a_{k,n} a_{k+1,n}^{-1}=1 \end{equation} for $j=n$ and $k\in
\mathbb Z$, since
$$
\begin{array}{l} x_n^k(yx_jy^{-1} x_j^{-1})x_n^{-k}
=\\
(x_n^kyx_n^{-(k+n)})(x_n^{k+n}x_jx_n^{-(k+n+1)}) (x_n^{k+n+1}
y^{-1} x_n^{-(k+1)})(x_n^{k+1}x_j^{-1}x_n^{-k})
 \end{array}
$$
and
$$ \begin{array}{l} x_n^k(y x_ny^{-1}x_n^{-1})x_n^{-k} = (x_n^k y
x_n^{-(k+n)}) (x_n^{k+n+1}y^{-1}x_n^{-(k+1)}).
\end{array}
$$

Similarly we can rewrite the relations $\bar R$ and denote the
result by $\tilde R$. Later  $\tilde R$ will be considered in more
detail and now let us show that $N$ is finitely generated. We have
$$N =
< a_{i,j}, i\in \mathbb Z, j\in \{2,\dots ,n\}\, \, |\, \,
(\ref{eq1}), (\ref{eq2}), \tilde R>. $$ Defining relations
(\ref{eq1}) and (\ref{eq2}) are equivalent to
\begin{equation} a_{0,n} = a_{i,n};\label{eq3}
\end{equation}
\begin{equation}
a_{0,n}a_{i+n,j}a_{0,n}^{-1} = a_{i,j}\label{eq4}.
\end{equation}
where $i\in \mathbb Z$ and $j\in \{2,\dots ,n\}$.

It follows from (\ref{eq4}) that $a_{i+kn,j} = a_{0,n}^{-k}
a_{i,j}a_{0,n}^{k}$ for $j\in \{2,\dots ,n-1\}, i\in\{0,\dots
,n-1\}$, and $k\in \mathbb Z$. Hence by Tietze transformations,
$$N = < a_{0,n}, a_{i,j}, i\in\{0,\dots ,n-1\}, j\in \{2,\dots ,n-1\} \,
\, |\, \, \hat R>,$$ where $\hat R$ is obtained from $\tilde R$ by
substitutions $ a_{k,n} = a_{0,n}$ and $a_{i+kn,j} =
a_{0,n}^{-k}a_{i,j}a_{0,n}^{k}$ in the words of $\tilde R$ for
$i\in\{0,\dots ,n-1\}, j\in \{2,\dots ,n-1\}$, and $k\in \mathbb
Z$.

So $N$ is generated by $a_{0,n}, a_{i,j}$ for $i\in\{0,\dots
,n-1\}$ and $j\in \{2,\dots ,n-1\}$ the number of which is
$n(n-1)+1 = (n-1)^2$. In particular, if $R$ is empty then $N$ is a
free group freely generated by these elements.

Now let us return to $\tilde R$. Each $r\in \bar R$ gives the set
of relations $\{ x_n^krx_n^{-k}, k\in \mathbb Z\}\subset \tilde
R$, where $\tilde r_k = x_n^k\cdot r\cdot x_n^{-k}$ for each $k\in
\mathbb Z$ is rewritten on the generators $a_{i,j}, i\in \mathbb
Z, j\in \{2,...,n\}$. Moreover it follows from the Reidemeister
rewriting process for $N$ that a word $\tilde r_l$ rewritten on
$\{ a_{i,j} \}$ coincides with $\tilde r_m$ after changing each
generator $a_{l+\nu, \mu}$ by $a_{m+\nu, \mu}$ in $\tilde r_l$,
since $a_{m+\nu,\mu} = x_n^{m-l}a_{l+\nu,\mu}x_n^{l-m}$ and
$\tilde r_m = x_n^{m-l}{\tilde r_l}x_n^{l-m}$. Then we have from
(\ref{eq3}) and (\ref{eq4}) that $\tilde
r_{j+kn}=a^{-k}_{0,n}{\tilde r_j}a^k_{0,n}$ for each $k\in \mathbb
Z$ and $j\in \{0,\dots ,n-1\}$. Therefore $\tilde r_{j+kn}$ is a
consequence of $\tilde r_j$ for each $r\in\bar R$, $k\in \mathbb
Z$ and $j\in\{0,\dots ,n-1\}$.  So
 for $k\in \mathbb Z\setminus \{0\}$ and
$j\in\{0,\dots ,n-1\}$ the relations $\tilde r_{j+kn}$ can be
removed from the set $\tilde R$ of defining relations of $N$.
Since $R$ is finite, the set $\{\tilde r_j = x_n^j\cdot r\cdot
x_n^{-j} | j\in\{0,\dots ,n-1\},
r\in R\}$ is finite. Hence $N$ is finitely presented, i.e., Theorem \ref{Main3} is proved. 

\section{Existence of non-Hopfian Hurwitz $C$-groups}
It is well-known (see, for example, \cite{L-S}) that the group
$$\tilde G = <a,t|t^{-1}{a^2}t = a^3>$$ is non-Hopfian. Therefore to prove
part (i) of Theorem \ref{Zar}, it is sufficient to show that
$\tilde G$ is an irreducible $C$-group. From the Tietze theorem we
have \begin{equation} \label{pres}
\begin{array}{ll} \tilde G
\simeq  & <a, t,x_1,x_2\, \, |\, \, x_1=t,x_2=ta,
{a^2}ta^{-2}=ta>\simeq \\ & <x_1, x_2\, \, |\, \,
(x_1^{-1}x_2)^2x_1(x_1^{-1}x_2)^{-2}=x_2>\simeq \\
 & <x_1,...,x_5 \, \, |\, \,  x_3 = x_1^{x_2}, x_3=
x_4^{x_1}, x_5 = x_4^{x_2}, x_5 = x_2^{x_1}>,
\end{array}\end{equation}
where $x_i^{x_j}= x_jx_ix_j^{-1}$, that is, $\tilde G$ is a
$C$-group. It is an irreducible $C$-group, since the
$C$-generators $x_1,\dots , x_5$ are conjugated to each other in
$\tilde G$.

To prove part (ii) of Theorem \ref{Zar}, consider the group
$$G=\tilde G\times <y\mid \emptyset >.$$ We have $G/G'\simeq \mathbb
Z^2$ and
$$
\begin{array}{ll} G\simeq & < x_1 , x_2, y\, |\,  yx_1 =
x_1y, x_2y=yx_2, (x_1^{-1}x_2)^2x_1(x_1^{-1}x_2)^{-2}=x_2
>\simeq \\ & <x_1,x_2,x_3,y \,\, \mid \,\, x_3=(x_1x_2)^{-1}y, yx_1 = x_1y, x_2y=yx_2,
\\ & (x_1^{-1}x_2)^2x_1(x_1^{-1}x_2)^{-2}=x_2>\simeq \\ & <x_1,x_2,x_3
\,\, \mid \,\, (x_1^{-1}x_2)^2x_1(x_1^{-1}x_2)^{-2}=x_2,
\\ & (x_1x_2x_3)x_i=x_i(x_1x_2x_3), \, i=1,2,3>.
\end{array}$$
Therefore $G$ is a Hurwitz $C$-group of degree three. Since the
group $\tilde G$ is non-Hopfian and $G = \tilde G\times<y\, \mid
\, \emptyset>$, the group $G$ is also non-Hopfian.

It is easy to see that for any $k\in \mathbb N$ the projective
Hurwitz groups $\bar G_{3,k}\simeq \tilde G\times \mathbb
Z/3k\mathbb Z$ and therefore they are also non-Hopfian groups.

As is known (see, for example, \cite{L-S}), non-Hopfian groups are
not residually finite.
\newline
{\bf Remark.} It follows from \cite{Ku1} that $\tilde G$ is also
the group of a 2-knot, since the graph of the last presentation in
(\ref{pres}) is a tree.

 \ifx\undefined\bysame
\newcommand{\bysame}{\leavevmode\hbox to3em{\hrulefill}\,}
\fi

\end{document}